\documentclass[oneside,a4paper]{amsart}
\usepackage{amssymb}
\usepackage{xspace}
\usepackage{subfig}
\usepackage{ifthen}
\usepackage[vmargin=35mm]{geometry}
\usepackage[pdfauthor={Dan Drake}, pdftitle={Bijections from weighted
  Dyck paths to Schroeder paths}]{hyperref}

\newtheorem{theorem}{Theorem}
\newtheorem{proposition}[theorem]{Proposition}
\newtheorem{lemma}[theorem]{Lemma}
\newtheorem{corollary}[theorem]{Corollary}
\newtheorem{conjecture}[theorem]{Conjecture}
\theoremstyle{definition}
\newtheorem*{definition}{Definition}

\newcommand{\twofone}[3]{{}_{2}F_1\left(\begin{matrix}#1&#2\\&#3\end{matrix}\right)}
\newcommand{\threeftwo}[6]{%
{}_{3}F_2\left(\begin{matrix}#1&#2&#3\\&#4&#5\end{matrix}\,;#6\right)}

\newcommand{\nonthmref}[2][Lemma]{\hyperref[#2]{#1~\ref*{#2}}}

\newcommand{\spthing}[1]{\textcolor{green}{#1}}

\newcommand{\N}{\mathbb{N}}

\newif\ifhavetikz

\usepackage{tikz}
\havetikztrue

\ifhavetikz
  \usetikzlibrary{decorations,decorations.pathmorphing,arrows}

  \pgfdeclarelayer{background}
  \pgfsetlayers{background,main}
  \tikzset{pathfragment/.style={scale=1.0,
  local bounding box=bbox,
  execute at end scope=
  {\begin{pgfonlayer}{background}
     \draw[help lines]
   (bbox.south west) grid (bbox.north east);
   \end{pgfonlayer}}}}

  \tikzset{pathstep/.style={thick}}

  \newcommand{\startpathat}[2]{
  \pgfmathsetmacro{\x}{#1}
  \pgfmathsetmacro{\y}{#2}}
  \newcommand{\up}{
  \filldraw (\x,\y) circle (2pt);
  \draw[pathstep] (\x,\y) -- ++(1,1);
  \pgfmathsetmacro{\x}{\x + 1} \pgfmathsetmacro{\y}{\y + 1}
  \filldraw (\x, \y) circle (2pt);}
  \newcommand{\down}{
  \filldraw (\x,\y) circle (2pt);
  \draw[pathstep] (\x,\y) -- ++(1,-1);
  \pgfmathsetmacro{\x}{\x + 1} \pgfmathsetmacro{\y}{\y - 1}
  \filldraw (\x, \y) circle (2pt);}
  \newcommand{\spdown}{
  \filldraw[green] (\x,\y) circle (2pt);
  \draw[pathstep,green] (\x,\y) -- ++(1,-1);
  \pgfmathsetmacro{\x}{\x + 1} \pgfmathsetmacro{\y}{\y - 1}
  \filldraw[green] (\x, \y) circle (2pt);}
  \newcommand{\horiz}{
  \filldraw (\x,\y) circle (2pt);
  \draw[pathstep] (\x,\y) -- ++(2,0);
  \pgfmathsetmacro{\x}{\x + 2}
  \filldraw (\x, \y) circle (2pt);}

  \tikzstyle{vertex}=[circle, draw=black, fill=black, inner sep=1]
  \newcommand{\drawedge}[3][]{\draw (#2) to [out=75, in=105,
    looseness=1.5, #1] (#3);}

  \pgfrealjobname{dyck-schroeder}
\fi

\newcommand{\inputtikz}[1]{%
\ifdefined\fakejobname%
  \ifthenelse{\equal{\jobname}{\detokenize{tikz/#1-out}}}{%
  \beginpgfgraphicnamed{tikz/#1-out}\input{tikz/#1.tex}\endpgfgraphicnamed}{}%
\else%
  \beginpgfgraphicnamed{tikz/#1-out}\input{tikz/#1.tex}\endpgfgraphicnamed%
\fi}

\newcommand{\Sch}{Schr\"oder\xspace}

\DeclareMathOperator{\pathlen}{pathlen}

\begin{document}

\title{Bijections from weighted Dyck paths to Schr\"oder paths}
\author{Dan Drake}
\address{Department of Mathematical Sciences\\Korea Advanced Institute
   of Science and Technology\\Daejeon, Korea}
\email{ddrake@member.ams.org}
\urladdr{http://mathsci.kaist.ac.kr/~drake}
\date{\today}

\begin{abstract}
  Kim and Drake used generating functions to prove that the number of
  $2$-distant noncrossing matchings, which are in bijection with little
  \Sch paths, is the same as the weight of Dyck paths in which downsteps
  from even height have weight $2$. This work presents bijections from
  those Dyck paths to little \Sch paths, and from a similar set of Dyck
  paths to big \Sch paths. We show the effect of these bijections on the
  corresponding matchings, find generating functions for two new classes
  of lattice paths, and demonstrate a relationship with $231$-avoiding
  permutations.
\end{abstract}

\maketitle

\section{Introduction and preliminaries}
\label{sec:intro}

This work begins with the work of Kim and the present author
~\cite{drake.kim:k-distant} in which they studied, among other things,
$2$-distant noncrossing matchings.
Such matchings---which will be defined shortly---are naturally
enumerated by little \Sch paths.
In the process of describing connections between $k$-distant noncrossing
matchings and orthogonal polynomials, Drake and Kim used generating
functions to show that little \Sch paths are equinumerous with a certain
set of labeled Dyck paths.
We present here a bijective proof of that fact; the bijection has a
number of interesting properties and is a consequence of a bijection
between big \Sch paths and a similar set of labeled Dyck paths.

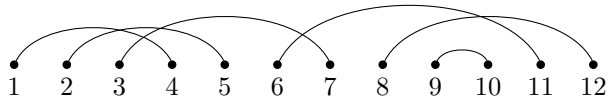
\begin{figure}[h]
  \centering
  \begin{tikzpicture}[scale=.7]
  \begin{scope}[yscale=.5]
    \foreach \a in {1,...,12}{
      \node (V\a) at (\a,0) [vertex, label=below:$\a$] {} ;
    }

    \drawedge{V1}{V4}
    \drawedge{V2}{V5}
    \drawedge{V3}{V7}
    \drawedge{V6}{V11}
    \drawedge{V8}{V12}
    \drawedge{V9}{V10}
  \end{scope}
\end{tikzpicture}

  \caption{A matching of $[12]$.}
  \label{fig:matching-example}
\end{figure}

We begin with definitions of the combinatorial objects mentioned above.
The notation $[n]$ refers to the set of positive integers from $1$ to
$n$.
A \textsl{matching of $[n]$} is a set of vertex-disjoint edges in the
complete graph on $n$ vertices so that every vertex is adjacent to
exactly one edge.
For our purposes, a matching can also be viewed as a permutation whose
cycles all have length $2$, or a set partition whose blocks all have
size $2$.
We will draw matchings by arranging the vertices horizontally and
drawing arcs, as in \autoref{fig:matching-example}.

Drake and Kim~\cite{drake.kim:k-distant} define a \textsl{$k$-distant
crossing} as a pair of arcs $(i_{1}, j_{1})$ and $(i_{2}, j_{2})$, with
$i_{1} < i_{2} < j_{1} < j_{2}$ and $j_{1} - i_{2} \ge k$.
The arcs $(6, 11)$ and $(8,12)$ of the matching in
\autoref{fig:matching-example} form a $3$-distant crossing; the arcs
$(1,4)$ and $(3,7)$ form a $1$-distant crossing.
A \textsl{$k$-distant noncrossing matching} is simply a matching with no
$k$-distant crossing.
The matching in \autoref{fig:matching-example} is $4$-distant
noncrossing.
(This notion of $k$-distant crossing is different from the $k$-crossings
of matchings studied by, for example, Chen et al.\
~\cite{chen.deng.ea:crossings}; their work concerns sets of $k$ mutually
crossing edges, and ignores the distance between vertices.)

\begin{figure}[h]
  \centering
  \begin{tikzpicture}[scale=.5]

  \begin{scope}
    \foreach \a in {1,...,14}{
      \node (V\a) at (\a,0) [vertex, label=below:$\a$] {} ;
    }

    \drawedge{V1}{V5}
    \drawedge{V2}{V3}
    \drawedge[looseness=.75]{V4}{V14}
    \drawedge{V6}{V8}
    \drawedge{V7}{V9}
    \drawedge{V10}{V13}
    \drawedge{V11}{V12}
  \end{scope}

  \begin{scope}[pathfragment, yshift=-5cm, xshift=0.5cm]
    \startpathat{0}{0}
    \up \up \down \horiz \up \horiz \down \up \up \down \down \down
  \end{scope}

  \begin{scope}[xshift=18cm]

    \begin{scope}
    \node at (0,0) [vertex] {};
    \draw (0,0) arc (180:90:1);

    \draw[<->] (0,-3/2) -- (0,-1/2);

    \startpathat{-1/2}{-3} \up
    \end{scope}

    \begin{scope}[xshift=3cm]
    \node at (0,0) [vertex] {};
    \draw (0,0) arc (0:90:1);

    \draw[<->] (0,-3/2) -- (0,-1/2);

    \startpathat{-1/2}{-2} \down
    \end{scope}

    \begin{scope}[xshift=6cm]
    \node at (-1/3,0) [vertex] {};
    \draw (-1/3,0) arc (180:90:1);
    \node at (1/3,0) [vertex] {};
    \draw (1/3,0) arc (0:90:1);

    \draw[<->] (0,-3/2) -- (0,-1/2);

    \startpathat{-1}{-5/2} \horiz
    \end{scope}

  \end{scope}

\end{tikzpicture}

  \caption{An example of the bijection between $2$-distant noncrossing
    matchings and little \Sch paths (left), and the correspondence
    between the edges incident to vertices of the matching and steps in
    the little \Sch path (right). }
  \label{fig:2dnc-matching-little-sch}
\end{figure}
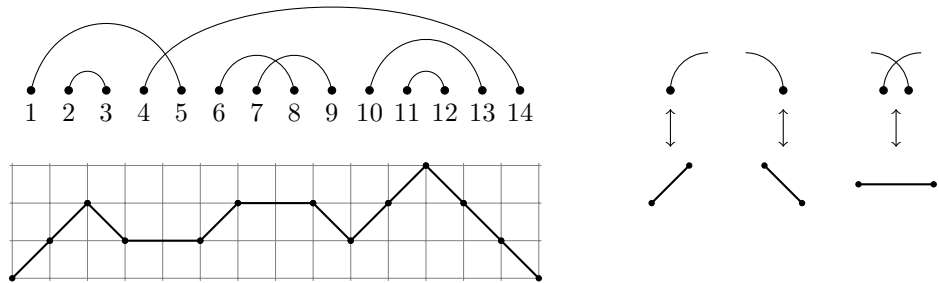

In a $2$-distant noncrossing matching, crossing edges are allowed as
long as the right vertex of the left edge is adjacent to the left
vertex of the right edge.
This fact allows us to describe a bijection from $2$-distant
noncrossing matchings to a certain class of lattice paths.
A lattice path of length $n$ is a sequence $(p_{0}, p_{1},\dots,p_{n})$
of points in $\N \times \N$; the $k$th step of the path is the pair
$(p_{k-1}, p_{k})$.
A step is called an \textsl{upstep} if the component-wise difference of
$p_{k} - p_{k-1}$ is $(1,1)$, and a \textsl{downstep} if the difference is
$(1, -1)$.
In this work, we will use paths with \textsl{double horizontal steps},
which is a pair of adjacent steps whose component-wise differences are
both $(1,0)$.
By a minor abuse of terminology, a double horizontal step will usually
be called a horizontal step.
A \textsl{little \Sch path} is a lattice path consisting of upsteps,
downsteps, and horizontal steps, such that no horizontal step occurs at
height zero.
See \autoref{fig:2dnc-matching-little-sch} for an example of such a path.

It is not difficult to describe a bijection from $2$-distant noncrossing
matchings to little \Sch paths: convert every vertex at the left end of
an arc to an upstep and every vertex at the right end of an arc to a
downstep---except for adjacent vertices involved in a crossing: convert
those vertices into a horizontal step.
This operation is a bijection because, given a little \Sch path, one can
recover the matching by drawing an opening half edge at every upstep,
two crossing half edges at every horizontal step, and a closing half
edge at every downstep.
Then connect every closing half-edge to the nearest opening half-edge to
create a matching.
See \autoref{fig:2dnc-matching-little-sch} for an example of this
correspondence.

The \textsl{little \Sch numbers} $s_{n}$ (sequence A1003 in the OEIS
~\cite{sloane:on-line}) count $2$-distant noncrossing matchings of
$[2n]$ and also little \Sch paths of length~$2n$.
If horizontal steps on the $x$-axis are allowed, one has a \textsl{big
\Sch path}; the number of such paths of length $2n$ is $S_{n}$, the
\textsl{big \Sch number} (sequence A6318) and it is well known that
$S_{n} = 2s_{n}$ for $n > 0$; see the next section and also Deutsch's
bijective proof~\cite{deutsch:bijective}.

We need several more definitions related to lattice paths.
The first is the step that \textsl{matches} a step.
For an upstep $u$, the matching step is the rightmost downstep to the
left of $u$ that leaves from the same height at which the upstep ends;
the definition for a downstep is similar.
For a horizontal step $h$ not on the $x$-axis, the matching step is the
leftmost downstep to the right of $h$ that leaves from the same height
as $h$; in the corresponding $2$-distant noncrossing matching, the
matching downstep corresponds to the rightmost vertex involved in the
two crossing edges.
For example, the step matching the first horizontal step in
\autoref{fig:2dnc-matching-little-sch} is the last downstep.
We will write paths using ``U'' for upsteps, ``HH'' for horizontal
steps, ``D'' for regular downsteps, and ``d'' for special downsteps,
which will be defined in \autoref{sec:descr-biject}.
The path in \autoref{fig:2dnc-matching-little-sch} is UUDHHUHHDUUDDD.

\subsection{Orthogonal polynomials and weighted Motzkin paths}
\label{sec:orth-polyn-weight}

When one has a sequence of positive numbers, in many cases it is
possible to describe that sequence as the moments of a sequence of
orthogonal polynomials.
In other words, given $\{\mu_{n}\}_{n \ge 0}$, define a measure (or a
linear functional on polynomials; the two are equivalent here) by
$\int x^{n} \,\mathrm{d} \mu = \mu_{n}$ and find polynomials $\{P_{n}(x)\}_{n
\ge 0}$ so that the integral
\begin{displaymath}
  \int P_{n}(x) P_{m}(x) \,\mathrm{d} \mu = 0
\end{displaymath}
when $n \ne m$ and is nonzero when $n = m$.
Many classical combinatorial sequences produce sequences of orthogonal
polynomials: the Catalan numbers produce Chebyshev polynomials of the
second kind, matching numbers produce Hermite polynomials, factorials
produce Laguerre polynomials, and so on.

Viennot described a completely combinatorial theory of orthogonal
polynomials ~\cite{viennot:theorie,viennot:combinatorial} in which the
moments of a sequence of orthogonal polynomials are expressed as
weighted Motzkin paths.
A Motzkin path is a lattice path that consists of upsteps, downsteps
and single horizontal steps (steps that move $(1,0)$); a weighted
Motzkin path has a weight $\lambda_ {n}$ associated with every downstep
leaving from height $n$, a weight $b_{n}$ for every horizontal step at
height $n$, and weight $1$ for all upsteps.
For many orthogonal polynomial moment sequences, the weights $b_{n}$
are zero, which means the corresponding moments may be described by
weighted Dyck paths; a \textsl{Dyck path} is just a Motzkin path with
no horizontal steps.

Drake and Kim~\cite{drake.kim:k-distant} showed that the number of
$2$-distant noncrossing matchings of $[2n]$---little \Sch numbers---is
the same as the total weight of weighted Dyck paths of length $2n$ in which
downsteps leaving from odd height have weight $1$, and downsteps leaving
from even height have weight $2$.
They proved this equality using equation (2) of Kim and Zeng
~\cite{kim.zeng:combinatorics} (or equation (1) of Vauchassade de
Chaumont and Viennot ~\cite{chaumont.viennot:polynomes}), which in the
present context is
\begin{equation}
  s_{n} = \sum_{k \ge 0} \frac{1}{n} \binom{n}{k} \binom{n}{k+1} 2^{k};
  \label{eq:little-sch-sum-of-narayana}
\end{equation}
in both works, the authors demonstrate that the sum above represents the
generating function for the weighted Dyck paths described above.
However, the sum also counts little \Sch paths, since
$\binom{n}{k}\binom{n}{k+1}/n$ is a \textsl{Narayana number} (sequence
A1263), which counts Dyck paths of length $2n$ with $k+1$ peaks and
$k$ ravines.
A peak is an upstep immediately followed by a downstep, and a ravine
is a downstep immediately followed by an upstep.
Between two consecutive peaks, there must be exactly one ravine, so
having $k+1$ peaks is equivalent to having $k$ ravines.
Any ravine can clearly be ``filled in'' and replaced with a horizontal
step, so a Dyck path with $k$ ravines corresponds to $2^{k}$ little
\Sch paths, which explains equation~\eqref{eq:little-sch-sum-of-narayana}.
On the other hand, any peak can be ``flattened'' into a horizontal
step, so we also have
\begin{equation}
  S_{n} = \sum_{k \ge 0} \frac{1}{n} \binom{n}{k} \binom{n}{k+1} 2^{k+1}
  \label{eq:big-sch-sum-of-narayana}
\end{equation}
because peaks can occur on the $x$-axis.
This provides one explanation for why there are twice as many big \Sch
paths as little ones.

\subsection{Plan of the paper}
\label{sec:plan-paper}

The aim of this work is to demonstrate a bijection from weighted Dyck
paths whose downsteps at even height have weight ~$2$ to little \Sch
paths.
That bijection will be a minor modification of a bijection from big \Sch
paths to a similar class of Dyck paths; both bijections will in turn be
consequences of more refined bijections between classes of little and
big hybrid paths, which are described in \autoref{sec:descr-biject}.
In \autoref{sec:hybrid-paths-matchings} we show the effect of those
bijections on the corresponding matchings, and then find generating
functions for little and big hybrid paths in
\autoref{sec:enum-hybr-paths}.
We finish by showing that our bijections are closely related to
$231$-avoiding permutations in \autoref{sec:hybrid-paths-231}.

\section{Description of the bijection}
\label{sec:descr-biject}

Instead of working with Dyck paths in which downsteps from even height
have weight $2$, we will work with Dyck paths in which such downsteps
may or may not be labeled ``special''; the two ideas are clearly
equivalent.
Such paths will be called \textsl{even-special Dyck paths} and abbreviated
``ESDPs''; \textsl{odd-special Dyck paths} (ODSPs) are defined similarly.

We will first describe a bijection $E_{\infty}$ from odd-special Dyck
paths to big \Sch paths---our desired bijection from even-special Dyck
paths to little \Sch paths will follow from a minor modification of
that bijection.
The bijection $E_{\infty}$ will be a consequence of a more refined
bijection $E$ between two classes of what we will call big hybrid
paths.
\textsl{Big hybrid paths} include odd-special Dyck paths and any path
obtained by applying $E$ to a big hybrid path.
To understand this recursive definition, we must define the map $E$.

\begin{definition}
  \label{def:e-bijection}
  Given a hybrid path, the map $E$ does nothing to the path if the path
  contains no special steps.
  Otherwise, given a hybrid path with $k$ horizontal steps, $E$ yields a
  hybrid path with $k+1$ horizontal steps by the following procedure.
  Find the leftmost special step in the hybrid path.
  If that step is preceded by an upstep, \textsl{flatten} the upstep and
  special downstep by replacing them with a horizontal step.
  If the special step is preceded by a downstep $d$, find the upstep $u$
  that matches $d$ and let $P$ be the (possibly empty) subpath between
  $u$ and $d$.
  Replace $u$ with a horizontal step, delete $d$, \textsl{slide} $P$ so
  that it follows the horizontal step, and make the original special step
  an ordinary downstep.
\end{definition}

\begin{figure}
  \centering
  \subfloat[When the special step is preceded by an upstep, $E$ and $e$
  simply ``flatten'' the two steps.]{
  \label{fig:e-bijection-flatten}
  \parbox{8.25cm}{\centering \begin{tikzpicture}[scale=.66]
\begin{scope}[pathfragment]
  \startpathat{0}{0} \up \spdown
\end{scope}

\draw[->] (2.5,0.5) -- (3.5,0.5);

\begin{scope}[pathfragment, xshift=4cm]
  \startpathat{0}{0} \horiz
  \draw (1, 1) circle (0);
\end{scope}
\end{tikzpicture}

}
}

  \subfloat[When the special step is preceded by a downstep, $E$ and $e$
  find the matching up step and ``slide'' $P$.]{
  \label{fig:e-bijection-slide}
  \begin{tikzpicture}[scale=.66]

\begin{scope}[pathfragment]
  \startpathat{0}{1} \up

  \draw[pathstep, decorate, decoration={snake}]
    (1,2) .. controls +(45:1) and +(135:1) .. +(3,0) node[midway,above,
    inner sep = 2mm]
    {$P$};

  \startpathat{4}{2} \down \spdown
\end{scope}

  \draw[->] (6.5,1.5) -- (7.5,1.5);

\begin{scope}[pathfragment, xshift=8cm]

  \startpathat{0}{1} \horiz
  \draw[pathstep, decorate, decoration={snake}]
    (2,1) .. controls +(45:1) and +(135:1) .. +(3,0) node[midway,above,
    inner sep=2mm]
    {$P$};
  \startpathat{5}{1} \down

  \draw (2.5, 3.5) circle (0);
\end{scope}

\end{tikzpicture}

}

\caption{The action of bijections $E$ and $e$ on the leftmost special
  down step in a hybrid path.}
  \label{fig:e-bijection}
\end{figure}
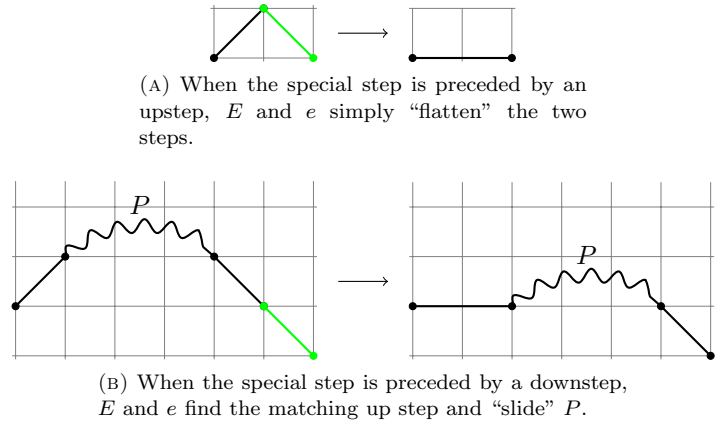

\autoref{fig:e-bijection} demonstrates the flatten and slide operations.
All the paths in \autoref{fig:E-forward-backward} are big hybrid paths.

The map $E$ clearly preserves the total number of special and horizontal
steps and, for paths with at least one special step, reduces the number
of special steps by one.
It is also a bijection:

\begin{theorem}
  The map $E$ is a bijection from the set of odd-special Dyck paths of
  length $n$ with no special steps to the set of big \Sch paths of
  length $n$ with no horizontal steps.
  It is also a bijection from the set of big hybrid paths of length
  $n$ with $j$ special steps and $k$ horizontal steps to the set of
  big hybrid paths with $j-1$ special steps and $k+1$ horizontal steps.
  \label{thm:E-is-a-bijection}
\end{theorem}

We will show that $E$ is a bijection by describing a procedure for
finding the horizontal step that was added last; the operation described
in \autoref{def:e-bijection} and \autoref{fig:e-bijection} is obviously
reversible if we know which horizontal step was added last.
Before giving the proof, let's see why this identification is not as
simple as it may sound.
The problem is that sometimes $E$ ``moves forward'' and sometimes $E$
``moves backward''.
\autoref{fig:E-forward-backward} shows what we mean by this.
A horizontal step may be created by $E$ to the left, to the right, or in
the middle of the existing horizontal steps, so the left- or rightmost
horizontal step need not be the last one added.

\begin{figure}[h]
  \centering
  \subfloat[The horizontal steps added by $E$ move ``backwards'' when
  doing repeated slide operations.]{
  \label{fig:E-slide-backward-ex}
  \begin{tikzpicture}[scale=.3]

\begin{scope}[pathfragment]
  \startpathat{0}{0} \up \up \up \up \up \spdown \down \spdown \down \spdown
\end{scope}

\begin{scope}[xshift=12cm]
  \draw[->] (-3/2, 2) -- (-1/2, 2);
  \node at (5,4) [above, font=\footnotesize] {$1$};
  \begin{scope}[pathfragment]
  \startpathat{0}{0} \up \up \up \up \horiz \down \spdown \down \spdown
  \end{scope}
\end{scope}

\begin{scope}[xshift=24cm]
  \draw[->] (-3/2, 2) -- (-1/2, 2);
  \node at (4,3) [above, font=\footnotesize] {$2$};
  \node at (6,3) [above, font=\footnotesize] {$1$};
  \begin{scope}[pathfragment]
  \startpathat{0}{0} \up \up \up \horiz \horiz \down \down \spdown
  \end{scope}
\end{scope}

\begin{scope}[xshift=36cm]
  \draw[->] (-3/2, 2) -- (-1/2, 2);
  \node at (2,1) [above, font=\footnotesize] {$3$};
  \node at (5,2) [above, font=\footnotesize] {$2$};
  \node at (7,2) [above, font=\footnotesize] {$1$};
  \begin{scope}[pathfragment]
  \startpathat{0}{0} \up \horiz \up \horiz \horiz \down \down
  \end{scope}
\end{scope}

\end{tikzpicture}

}

  \subfloat[The horizontal steps added by $E$ move ``forwards'' when
  doing repeated flatten operations.]{
  \label{fig:E-flatten-forward-ex}
  \begin{tikzpicture}[scale=.3]

  \begin{scope}[pathfragment]
  \startpathat{0}{0} \up \up \up \spdown \up \spdown \up \spdown \down \down
  \end{scope}

\begin{scope}[xshift=12cm]
  \draw[->] (-3/2, 2) -- (-1/2, 2);
  \node at (3,2) [above, font=\footnotesize] {$1$};
  \begin{scope}[pathfragment]
  \startpathat{0}{0} \up \up \horiz \up \spdown \up \spdown \down \down
  \end{scope}
\end{scope}

\begin{scope}[xshift=24cm]
  \draw[->] (-3/2, 2) -- (-1/2, 2);
  \node at (3,2) [above, font=\footnotesize] {$1$};
  \node at (5,2) [above, font=\footnotesize] {$2$};
  \begin{scope}[pathfragment]
  \startpathat{0}{0} \up \up \horiz \horiz \up \spdown \down \down
  \end{scope}
\end{scope}

\begin{scope}[xshift=36cm]
  \draw[->] (-3/2, 2) -- (-1/2, 2);
  \node at (3,2) [above, font=\footnotesize] {$1$};
  \node at (5,2) [above, font=\footnotesize] {$2$};
  \node at (7,2) [above, font=\footnotesize] {$3$};
  \begin{scope}[pathfragment]
  \startpathat{0}{0} \up \up \horiz \horiz \horiz \down \down
  \end{scope}
\end{scope}

\end{tikzpicture}

}

  \subfloat[Horizontal steps can also be added between existing horizontal
  steps.]{
  \label{fig:E-slide-middle-ex}
  \begin{tikzpicture}[scale=.3]

  \begin{scope}[pathfragment]
  \startpathat{0}{0} \up \up \up \spdown \down \up \up \spdown \down \spdown
  \end{scope}

\begin{scope}[xshift=12cm]
  \draw[->] (-3/2, 2) -- (-1/2, 2);
  \node at (3,2) [above, font=\footnotesize] {$1$};
  \begin{scope}[pathfragment]
  \startpathat{0}{0}  \up \up \horiz \down \up \up \spdown \down \spdown
  \end{scope}
\end{scope}

\begin{scope}[xshift=24cm]
  \draw[->] (-3/2, 2) -- (-1/2, 2);
  \node at (3,2) [above, font=\footnotesize] {$1$};
  \node at (7,2) [above, font=\footnotesize] {$2$};
  \begin{scope}[pathfragment]
  \startpathat{0}{0} \up \up \horiz \down \up \horiz \down \spdown
  \end{scope}
\end{scope}

\begin{scope}[xshift=36cm]
  \draw[->] (-3/2, 2) -- (-1/2, 2);
  \node at (3,2) [above, font=\footnotesize] {$1$};
  \node at (8,1) [above, font=\footnotesize] {$2$};
  \node at (6,1) [above, font=\footnotesize] {$3$};
  \begin{scope}[pathfragment]
  \startpathat{0}{0} \up \up \horiz \down \horiz \horiz \down
  \end{scope}
\end{scope}

\end{tikzpicture}

}

  \caption{The horizontal steps created by $E$ are not necessarily added
    left-to-right. The three paths on the far right look similar, but
    their horizontal steps were added in different orders.}
  \label{fig:E-forward-backward}
\end{figure}

One may think that, since horizontal steps from slides are always
created at odd height and horizontal steps from flattenings are created
at even height, it might be possible to use that information to identify
the last-added step, but since slides change the height of parts of the
path by one, a simple examination of odd and even heights will not
suffice.

\begin{proof}[Proof of \autoref{thm:E-is-a-bijection}]
  The first statement of the theorem is trivial, as it is saying that
  $E$ acts as the identity on the set of Dyck paths.
  For the second statement, we must show that it is possible to identify
  which horizontal step was added last.
  This can be done with the following procedure.

  Partition the path into subpaths that consist of either a sequence
  of non-horizontal steps, or a horizontal step, its matching step, and
  all steps in between.
  The only part of a path altered by $E$ when adding a horizontal step
  is between the horizontal step and its matching downstep, so if a
  horizontal step~$b$ is to the right of the downstep matching a
  horizontal step~$a$, then $b$ must have been added after ~$a$.
  (This is a special case of \nonthmref{thm:21-order-adding-sp-steps}.)
  This fact tells us that the last-added horizontal step must be in the
  rightmost such subpath that contains a horizontal step.
  Call that subpath the \textsl{first active subpath}.
  \autoref{fig:partition-by-hh} illustrates this partitioning process.

  \begin{figure}[ht]
    \centering
    \newcommand{\pathpartition}[2]{
\draw (#1,-1/2) -- (#1,7/2);
\draw (#2,-1/2) -- (#2,7/2);}

\begin{tikzpicture}[scale=.45]

\begin{scope}[pathfragment]

\startpathat{2}{0}

\up \down

\horiz \horiz \up

\up

\horiz \up \down \down

\horiz \down

\up

\horiz \up \horiz \up \down \down \down

\up \up \down \spdown
\end{scope}

\pathpartition{4}{9}
\pathpartition{10}{15}
\pathpartition{15}{18}
\pathpartition{19}{28}

\end{tikzpicture}

    \caption{The partitioning process to find the first active subpath,
      which is the rightmost subpath with a horizontal step.}
    \label{fig:partition-by-hh}
  \end{figure}

  If the first active subpath starts with a horizontal step on the
  $x$-axis, then, because horizontal steps on the $x$-axis can only be
  created with a flatten operation, the rightmost horizontal step in the
  subpath must be the last-added step.

  Otherwise, we may assume the first active subpath starts with a
  horizontal step at some positive height.
  If that step is at odd height, that step is the last-added
  horizontal step, because the slide operation of $E$ creates horizontal
  steps at odd height, and, as seen in
  \autoref{fig:E-slide-backward-ex}, as one moves forward along a
  sequence of downsteps, some of which are special steps, $E$ creates
  horizontal steps at the beginning of the first active subpath.

  If the step at the beginning of the first active subpath is at even
  height, we must partition the path again.
  Now partition the first active subpath into sequences of horizontal
  steps at the same height as the original horizontal step and subpaths
  that begin with an upstep and end at the downstep matching the upstep.
  Call these two kinds of sequences \textsl{valleys} and \textsl{hills},
  respectively.
  Using the same reasoning as before, the last-added horizontal step
  must be in the rightmost hill or valley that contains a horizontal step.
  Call that hill or valley the \textsl{second active subpath}.
  In the subpath of \autoref{fig:hills-valleys}, the final hill is the
  second active subpath.

  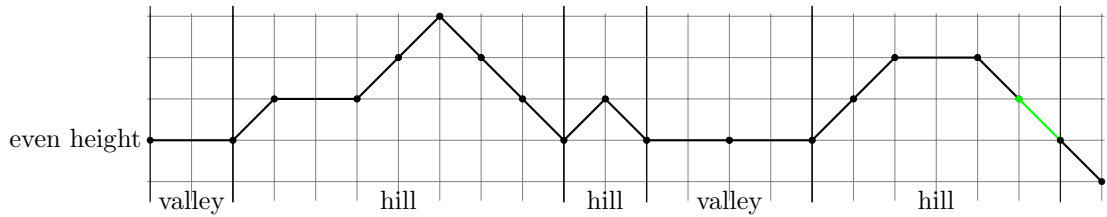
\begin{figure}[h]
    \centering
    \newcommand{\pathpartition}[3]{
\draw (#1,1/2) -- (#1,5.25);
\path (#1,1) -- (#2,1) node[midway, below] {#3};
\draw (#2,1/2) -- (#2,5.25);}

\begin{tikzpicture}[scale=.55]

\node at (0,2) [left] {even height};

\begin{scope}[pathfragment]

\startpathat{0}{2}

\horiz

\up \horiz \up \up \down \down \down

\up \down

\horiz \horiz

\up \up \horiz \down \spdown

\down

\node at (1/2, 3/4) {}; 
\end{scope}

\pathpartition{0}{2}{valley}
\pathpartition{2}{10}{hill}
\pathpartition{10}{12}{hill}
\pathpartition{12}{16}{valley}
\pathpartition{16}{22}{hill}

\end{tikzpicture}

    \caption{Partitioning the first active subpath into hills and
      valleys. The second active subpath is the rightmost hill or valley
      with a horizontal step.}
    \label{fig:hills-valleys}
  \end{figure}

  If the second active subpath is a valley, the rightmost step in the
  valley is the most recently added horizontal step because steps in a
  valley must come from the flattening operation.

  If the second active subpath is a hill, we recursively use the
  procedure described here to identify the last-added step within that
  hill.
  Since the hill begins at even height, the path is of the same form as
  the hybrid paths we began with.

  Since the paths have finite length and the recursion step uses a
  shorter path than it started with, this procedure always finishes, and
  since the ``exit points'' always identify what must be the most
  recently added step, the procedure as a whole will identify the
  last-added horizontal step of a hybrid path.
\end{proof}

For example, with the second active subpath in
\autoref{fig:hills-valleys}, we would recursively use the procedure on
the final hill.
The procedure in the proof above, given that hill (UUHHDd) as a single
path would identify HHD as the first active subpath then, since the
horizontal step is at even height, partition again and identify the
valley HH as the second active subpath, and finally declare that single
horizontal step as the most recently added horizontal step.
In \autoref{fig:e-bijection-ex}, step $11$ is the last-added
horizontal step.

\subsection{Consequences of the bijection}
\label{sec:consequences-of-E}

If we start with an odd-special Dyck path, we can use $E$ to iteratively
``evolve'' the path into a big \Sch path.
(In fact, we use $E$ to suggest the word ``evolve''.)
Let $E_{\infty}$ be the resulting map from odd-special Dyck paths to big
\Sch paths.
Since $E$ is a bijection and preserves the total number of special
steps and horizontal steps, we have the following corollary of
\autoref{thm:E-is-a-bijection}.

\begin{corollary} The map $E_{\infty}$ is a bijection from odd-special
Dyck paths of length $n$ with $k$ special steps to big \Sch paths of
length $n$ with $k$ horizontal steps.
  \label{thm:E-infinity-is-a-bijection}
\end{corollary}

The operation described in \autoref{def:e-bijection} and
\autoref{fig:e-bijection} does not refer to the parity of the heights
of the special steps, so we may use it with even-special Dyck paths.
Define the map $e$ the same way as $E$, but starting with even-special
Dyck paths.
\textsl{Little hybrid paths} are defined analogously to big hybrid
paths.
By simply switching ``odd'' and ``even'' in the proof of
\autoref{thm:E-is-a-bijection} and ignoring the possibility of
sequences of horizontal steps on the $x$-axis, we see that $e$ is also
a bijection:

\begin{corollary}
  The map $e$ is a bijection from the set of even-special Dyck paths of
  length $n$ with no special steps to the set of little \Sch paths of
  length $n$ with no horizontal steps.
  It is also a bijection from the set of little hybrid paths of length $n$ with
  $j$ special steps and $k$ horizontal steps to the set of hybrid paths
  with $j-1$ special steps and $k+1$ horizontal steps.
  \label{thm:e-is-a-bijection}
\end{corollary}

By defining $e_{\infty}$ analogously to $E_{\infty}$, we accomplish
our goal of showing bijectively that the little \Sch numbers enumerate
even-special Dyck paths of length $n$:

\begin{corollary}
  The map $e_{\infty}$ is a bijection from even-special Dyck paths of
  length $n$ with $k$ special steps to little \Sch paths of length $n$
  with $k$ horizontal steps.
  \label{thm:e-infinity-is-a-bijection}
\end{corollary}

\begin{figure}[h]
  \centering
  \begin{tikzpicture}[scale=.3]

\begin{scope}[pathfragment]
  \startpathat{0}{0} \up \up \up \up \spdown \up \spdown \down \up
  \down \down \up \spdown 

  \up \up \up \up \up \down \up \spdown \up \up \up \spdown \down
  \spdown \down \spdown \down \spdown \down
\end{scope}

\draw[->] (16,-1) -- node[right] {$e_{\infty}$} (16,-3);

\begin{scope}[yshift=-8cm]
  \begin{scope}[pathfragment]
    \startpathat{0}{0} \up \up \up \horiz \horiz \down \up \down \down
    \horiz \up \horiz \up \horiz \up \down \horiz \up \horiz \horiz
    \down \down \down \down
  \end{scope}
\end{scope}

\end{tikzpicture}

  \caption{An example of the bijection $e_{\infty}$.}
  \label{fig:e-bijection-ex}
\end{figure}
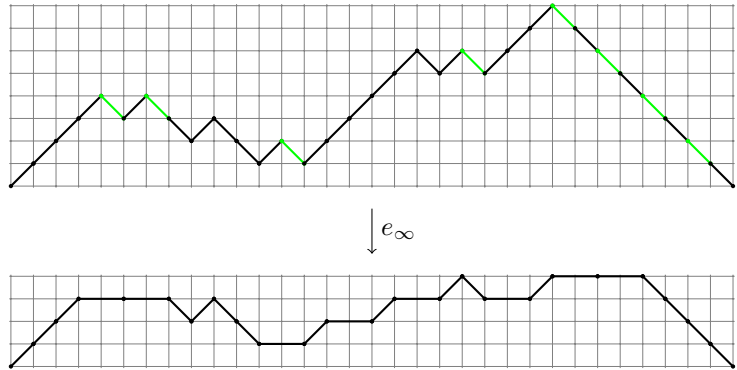

\autoref{fig:e-bijection-ex} shows an example of $e_{\infty}$.

Using the reasoning behind equations
\eqref{eq:little-sch-sum-of-narayana} and
\eqref{eq:big-sch-sum-of-narayana}, which counted \Sch paths by changing
peaks or ravines in Dyck paths into horizontal steps, the bijections
above imply that the number of even-special Dyck paths of length $2n$
with $j$ special steps is, for positive $n$,
\begin{equation}
  \sum_{k \ge 0} N(n,k) \binom{k}{j}
  =
  N(n, j)\,\twofone{j-n}{j-n+1}{j+2}
  =
  \frac{1}{n}\binom{n}{j}\binom{2n-j}{n+1},
  \label{eq:esdps-with-j-sp}
\end{equation}
where $N(n, k)$ is again a Narayana number and the ${}_{2} F_{1}$
notation is a hypergeometric function evaluated at one, which we can sum
with the Chu-Vandermonde identity.
The above triangle of numbers is sequence A126216.
Similarly, the number of odd-special Dyck paths of length $2n$ with
$j$ special steps is, for positive $n$,
\begin{equation}
  \sum_{k \ge 0} N(n,k) \binom{k+1}{j}
  =
  N(n, j-1)\,\twofone{j-n}{j-n-1}{j}
  =
  \frac{1}{n-j+1}\binom{n}{j}\binom{2n-j}{n};
  \label{eq:osdps-with-j-sp}
\end{equation}
the middle expression is not defined when $j=0$, but in that case, the
sum on the left is just the sum of the Narayana numbers---a Catalan
number---so the rightmost expression is correct for all nonnegative $j$.
The triangle in equation \eqref{eq:osdps-with-j-sp} is sequence A60693.

\section{Hybrid paths as matchings}
\label{sec:hybrid-paths-matchings}

This work began with an investigation of certain matchings, and since
little hybrid paths were developed to describe our bijection, it is
fitting that we examine the connection between little hybrid paths and
matchings.
We already know that little \Sch paths correspond to $2$-distant
noncrossing matchings, so first we will describe an interpretation
of even-special Dyck paths.
Using the bijection from \Sch paths to matchings, a Dyck path with no
special steps corresponds to a noncrossing matching, so it is reasonable
to interpret special steps in the path as special edges in the matching.
For example, the path UUdUUUDDdD corresponds to the noncrossing
matching $\{(1,10), \spthing{(2,3)}, \spthing{(4,9)}, (5,8), (6,7)\}$
in which the edges between $2$ and $3$ and between $4$ and $9$ are
special.

To interpret paths with both special steps and horizontal steps and
understand the action of $e$ in terms of paths, we need to define
nesting.
An edge $(a,b)$ in a matching \textsl{nests} the edge $(c,d)$ if $a <
c < d < b$.
An edge $a$ in a matching \textsl{immediately nests} edge $b$ if $a$
nests $b$, and any other edge that nests $b$ also nests $a$.

Before describing the action of the flatten and slide operations on
``hybrid matchings'', we need one observation.

\begin{lemma}
  Let $h$ be a double horizontal step in a little hybrid path.
  Let $d$ be the downstep matching $h$ and $u$ the upstep matching $d$.
  The step in the path corresponding to the rightmost (respectively,
  leftmost) vertex involved in the $1$-distant crossing at $h$ is either
  $d$ (resp., $u$) or the leftmost (resp., rightmost) horizontal step to
  the right (resp., left) of $h$ which is at the same height as $h$,
  whichever is closer to $h$.
  \label{thm:horiz-to-crossing-matching-steps}
\end{lemma}

\begin{proof}
  A key idea in this proof is that a sequence of steps in a hybrid
  path that begins with an upstep and ending with the matching
  downstep---what we called a hill in the proof of
  \autoref{thm:E-is-a-bijection}--- corresponds to a set of vertices
  in the matching that form a ``submatching''.
  Let $h$ be a double horizontal step in a little hybrid path, and
  partition the path as we did to find the second active subpath in the
  proof of \autoref{thm:E-is-a-bijection} (but ignore the height of
  $h$).
  The first step in $h$ corresponds to an opening half edge $x$.
  To what vertex will $x$ be connected?
  Any hill to the right of $h$ corresponds to a group of vertices that
  form a submatching, and hence $x$ will not be connected to any of
  them.
  If there is a double horizontal step between $h$ and $d$, then $x$
  will be connected to the vertex corresponding to the second step of the
  leftmost such double horizontal step; otherwise, $x$ will be connected
  to the vertex corresponding to $d$.
  The proof for the ``respectively'' part of the statement is similar.
\end{proof}

For example, in \autoref{fig:2dnc-matching-little-sch}, the horizontal
step at positions $4$ and $5$ corresponds to the $1$-distant crossing in
the matching in the same position; step $14$ in the \Sch path is the
downstep that matches the horizontal edge, so the right vertex of the
edge incident to vertex $4$ is vertex $14$; and step $1$ in the path is
the upstep that matches step $14$, so vertex $1$ is the left vertex of
the edge incident to vertex $5$.
Another example is UHHUDHHD; vertex $2$ in the corresponding matching
(which is $\{(1, 3), (2, 7), (4, 5), (6, 8)\}$) connects to vertex $7$
because the second horizontal step is at the same height as the first,
and is between the first horizontal step and its matching downstep.

Now we can describe the corresponding action of the bijection $e$ on
hybrid matchings.

\begin{figure}[ht]
  \centering
  \subfloat[The flatten operation when applied to a matching just ``swaps
  tails''.]{
  \label{fig:e-bijection-matching-flatten}
  \begin{tikzpicture}[scale=.5]
  \begin{scope}
    \node (Va) at (0,0) [vertex] {};
    \node (Vb) at (10,0) [vertex] {};
    \node (Vc) at (4,0) [vertex] {};
    \node (Vd) at (5,0) [vertex] {};

    \drawedge[looseness=1]{Va}{Vb}
    \drawedge[green]{Vc}{Vd}

    \draw[pathstep, decorate, decoration={snake}]
     (.75,0) .. controls +(20:1) and +(160:1) .. (3.25,0)
     node[midway,above, inner sep = 2mm] {$M_{1}$};

    \draw[pathstep, decorate, decoration={snake}]
     (5.75,0) .. controls +(20:1) and +(160:1) .. (9.25,0)
     node[midway,above, inner sep = 2mm] {$M_{2}$};
  \end{scope}

  \draw[->] (11,1) -- (13,1);

  \begin{scope}[xshift=14cm]
    \node (Va) at (0,0) [vertex] {};
    \node (Vb) at (10,0) [vertex] {};
    \node (Vc) at (4,0) [vertex] {};
    \node (Vd) at (5,0) [vertex] {};

    \drawedge{Va}{Vd}
    \drawedge{Vc}{Vb}

    \draw[pathstep, decorate, decoration={snake}]
     (.75,0) .. controls +(20:1) and +(160:1) .. (3.25,0)
     node[midway,above, inner sep = 2mm] {$M_{1}$};

    \draw[pathstep, decorate, decoration={snake}]
     (5.75,0) .. controls +(20:1) and +(160:1) .. (9.25,0)
     node[midway,above, inner sep = 2mm] {$M_{2}$};
  \end{scope}
\end{tikzpicture}

}

  \subfloat[The slide operation when applied to a matching. Note that the
  edge from $c$ to $b-1$ may actually be several edges, as in UUUHHDdD,
  where $b-1$ is $6$ and $c$ is $3$.]{
  \label{fig:e-bijection-matching-slide}
  \newcommand{\sameinboth}{
  \begin{scope}
    \node (Va) at (0,0) [vertex, label=below:$a\vphantom{1}$] {};
    \node (Vb) at (10,0) [vertex, label=below:$b$] {};
    \node (Vc) at (5,0) [vertex, label=below:$c\vphantom{1}$] {};
  \end{scope}

    \draw[pathstep, decorate, decoration={snake}]
     (.75,0) .. controls +(20:1) and +(160:1) .. (4.25,0)
     node[midway,above, inner sep = 2mm] {$M_{1}$};
}

\newcommand{\Mtwo}{
    \draw[pathstep, decorate, decoration={snake}]
     (5.75,0) .. controls +(20:1) and +(160:1) .. (8.25,0)
     node[midway,above, inner sep = 2mm] {$M_{2}$};
}

\begin{tikzpicture}[scale=.55, every label/.style={font=\footnotesize}]
  \begin{scope}
    \sameinboth
    \node (Vd) at (9,0) [vertex, label=below:$b-1$] {};

    \drawedge[looseness=1,green]{Va}{Vb}
    \drawedge{Vc}{Vd}

    \node at (6,0) [vertex, draw=white, fill=white, label=below:$c+1$] {};

    \Mtwo
  \end{scope}

  \draw[->] (11,1) -- (13,1);

  \begin{scope}[xshift=14cm]
    \sameinboth
    \node (Vd) at (6,0) [vertex, label=below:$c+1$] {};

    \drawedge{Va}{Vd}
    \drawedge{Vc}{Vb}

    \begin{scope}[xshift=1cm] \Mtwo \end{scope}
  \end{scope}
\end{tikzpicture}

}

  \caption{The effect of $e$ on matchings. This is the matchings version
  of \autoref{fig:e-bijection}.}
  \label{fig:e-bijection-matchings}
\end{figure}
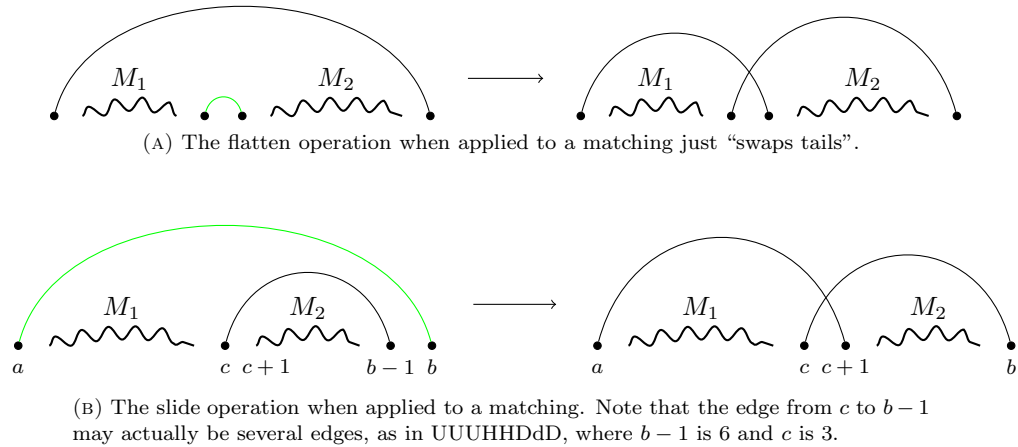

\begin{theorem}
  The analogue of the flatten operation for hybrid paths works as
  follows on matchings: given a special edge connecting vertices $c$ and
  $c+1$, find the edge that immediately nests that special edge; say it
  connects $a$ and $b$.
  Then ``swap the tails'': replace the edges $(c,c+1)$ and $(a,b)$ with
  edges $(c,b)$ and $(a,c+1)$.
  \label{thm:e-matchings-flatten}
\end{theorem}

\begin{proof}
  The special downstep in the path (which is immediately preceded by an
  upstep, since we are doing a flatten operation) becomes a horizontal
  edge.
  So there will be $1$-distant crossing at vertices $c$ and $c+1$.
  We only need to find the other two vertices involved; because of the
  way matchings are constructed from paths and
  \nonthmref{thm:horiz-to-crossing-matching-steps}, those other two
  vertices are $a$ and $b$.
\end{proof}

To define the analogue of the slide operation, we need to define the
\textsl{transitive left endpoint} of an edge.
Given an edge $e$, the transitive left endpoint of that edge is simply
the left endpoint of $e$---unless the edge is the right edge in a
$1$-distant crossing; then the transitive left endpoint of $e$ is the
transitive left endpoint of the left edge in the crossing.
For example, the transitive left endpoint of the edge $(6,8)$ in the
matching $\{(1,5), (2,3), (4,7), (6,8)\}$ is $1$. In little hybrid
paths, the transitive left endpoint corresponds to finding the upstep
that matches a downstep; because of
\nonthmref{thm:horiz-to-crossing-matching-steps}, a matching
upstep-downstep pair may not correspond to the left and right vertices
of a single edge.

\begin{theorem}

  The analogue of the slide operation for paths works as follows on
  matchings: given a special edge $(a,b)$, let $c$ be the transitive left
  vertex of the edge incident to vertex $b-1$.
  The slide operation on hybrid paths corresponds to replacing $(a,b)$
  and $(c,b-1)$ with edges $(a,c+1)$ and $(c,b)$ and sliding all half edges
  incident to vertices from $c+1$ to $b-2$ to the right by one
  vertex.
  \label{thm:e-matchings-slide}
\end{theorem}

\begin{proof}
  Since we are doing a slide operation, the special step in the path
  must be preceded by a ordinary downstep, which means $b-1$ must be the
  right vertex of an ordinary edge nested by the special edge.
  The new $1$-distant crossing created by the slide operation will be at
  the upstep that matches the downstep at $b-1$, which as we saw above
  is the transitive left endpoint of the edge incident to $b-1$ in the
  matching. We create the new crossing at vertices $c$ and $c+1$; the
  slide operation on the path moves all steps from $c+1$ to $b-2$ to
  $c+2$ to $b-1$, so all half edges incident to vertices from $c+1$ to
  $b-2$ are moved to the right one vertex.
\end{proof}

\autoref{fig:e-bijection-matchings} demonstrates these two operations
for matchings.

\section{Enumeration of hybrid paths}
\label{sec:enum-hybr-paths}

Having defined and used hybrid paths it is natural to wonder just many
of them there are.
All big \Sch paths and odd-special Dyck paths are big hybrid paths, and
Dyck paths, which are counted by the Catalan number $C_{n}$, are both
big \Sch paths and OSDPs, so there are certainly at least $2 S_{n} -
C_{n}$ big hybrid paths, but there are paths such as HHUd which are
neither \Sch paths nor odd-special Dyck paths.
\autoref{tab:hybrid-paths} shows the number of all hybrid paths,
little and big, for some small values of $n$.

\begin{table}[ht]
  \centering
  \begin{tabular}{lrrrrrrrrrrrr}
   $n$:  & 0   &  2  &   4  &   6  &    8  &    10  &    12  &     14  &
   16  &      18  &       20  & 22\\
  \hline
   little:  & 1 &  1  &   4  &  18  &   87  &   439  &  2278  &  12052  &
   64669  &  350733  &  1918152  &  10560678\\
   big:  & 1   &  3  &  11  &  47  &  219  &  1075  &  5459  &  28383  &
   150131  &  804515  &  4355163  & 23768079\\
  \end{tabular}

  \caption{The number of little and big hybrid paths.}
  \label{tab:hybrid-paths}
\end{table}

One way to count hybrid paths is to begin with even- and odd-special
Dyck paths with $j$ special steps, which are counted in equations
\eqref{eq:esdps-with-j-sp} and \eqref{eq:osdps-with-j-sp}; repeatedly
applying $e$ or $E$ to such a Dyck path will produce $j$ hybrid paths.
Multiplying those equations by $j+1$ and summing over $j$ yields, for
little hybrid paths,
\begin{equation}
  C_{n}\,\threeftwo{-n}{-n+1}{2}{-2n}{1}{-1},
  \label{eq:little-hyb-3f2}
\end{equation}
where $C_{n}$ is a Catalan number and the hypergeometric function is now
evaluated at $-1$. Similarly, the number of big hybrid paths is
\begin{equation}
  C_{n}\,\threeftwo{-n}{-n-1}{2}{-2n}{1}{-1}.
  \label{eq:big-hyb-3f2}
\end{equation}

Another way to enumerate these paths is to find their generating
functions.
Let $L(x)$ and $B(x)$ be the generating functions for little and big
hybrid paths, respectively.
We will use the following generating functions: $E(x)$ and $O(x)$ for
even- and odd-special Dyck paths, and $s(x)$ and $S(x)$ for little and
big \Sch paths.
Of course, we already know that $E(x) = s(x)$, $O(x) = S(x)$, $S(x) = 2
s(x) - 1$, and
\begin{displaymath}
  s(x) = \frac{2}{1 + x + \sqrt{x^2 - 6x + 1}},
\end{displaymath}
but it will be helpful to use different names to keep different types of
paths separate.
In all the generating functions considered here, paths of length $2n$
are weighted by $x^{n}$.

\begin{theorem}
Let $R = \sqrt{x^{2} - 6x + 1}$. The ordinary generating function
for little hybrid paths is
\begin{equation}
  L(x) = \frac{R+1-x}{2} \cdot
         \frac{2(R+x)}{R (R+x+1)} \cdot
         \frac{2}{R+x+1}
  \label{eq:little-hybrid-gf-R}
\end{equation}
and the ordinary generating function for big hybrid paths is
\begin{equation}
  B(x) = \left(\frac{R+1-x}{2} - 1 + x -
               \frac{R(R+x)}{2} + \frac{3}{2} - \frac{3x}{2} \right)
         \frac{2(R+x)}{R (R+x+1)} \cdot
         \frac{2}{R+x+1}.
\label{eq:big-hybrid-gf-R}
\end{equation}
More explicitly, we have
\begin{equation}
  L(x) = \frac{1 -5x + \sqrt{x^{2} -6x +1}}
{\left(
(x-1)^{2} + (x+1)\sqrt{x^{2} -6x+1}
\right) \sqrt{x^{2}-6x+1}}
\label{eq:little-hybrid-gf}
\end{equation}
and
\begin{equation}
  B(x) = \frac{7x - 2x^{2} - 1 + \sqrt{x^{2} - 6x + 1}}{2x \sqrt{x^{2}
      - 6x + 1}} - 1.
\label{eq:big-hybrid-gf}
\end{equation}
\end{theorem}

\begin{proof}
  We will decompose little and big hybrid paths to express $L$ and $B$
  in terms of each other and then solve the system.
  A key idea is that raising a big hybrid path up one unit and
  sandwiching it between and upstep and downstep yields a valid little
  hybrid path, since all the step height parities have effectively been
  reversed.
  Doing the same thing to a little hybrid path yields a big hybrid path,
  although the resulting path will never have a horizontal step at height
  ~$1$.

  Every nonempty little hybrid path may be decomposed into an upstep, a
  big hybrid path $P$, a downstep, then a little hybrid path $Q$, as shown
  in \autoref{fig:little-hybrid-decomp}.
  Any pair $P$ and $Q$ is allowed, unless $P$ has a special step and $Q$
  has a horizontal step.

  \begin{figure}[h]
    \centering
    \begin{tikzpicture}[scale=3/5]

\begin{scope}[pathfragment]

\startpathat{0}{0} \up

\draw[pathstep, decorate, decoration={snake}]
    (1,1) .. controls +(45:1) and +(135:1) .. +(4,0) node[midway,above,
    inner sep = 2mm]
    {big hybrid $P$};

\startpathat{5}{1} \down

\draw[pathstep, decorate, decoration={snake}]
    (6,0) .. controls +(45:1) and +(135:1) .. +(4,0) node[midway,above,
    inner sep = 2mm]
    {little hybrid $Q$};

\node[vertex] at (10,0) {};

\end{scope}

\end{tikzpicture}

    \caption{A decomposition of a nonempty little hybrid path.}
    \label{fig:little-hybrid-decomp}
  \end{figure}
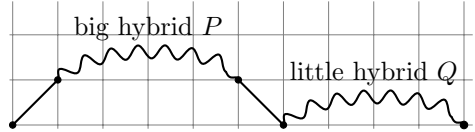

  Assume that $P$ has a special step, so that $Q$ has no horizontal
  step.
  Since $S(x)$ counts big hybrid paths with only upsteps, downsteps, and
  horizontal steps---in other words, with no special steps---the
  generating function for big hybrid paths with a special step is $B(x) -
  S(x)$; similar reasoning shows that the generating function for little
  hybrid paths with no horizontal steps is simply $E(x)$.
  Thus the generating function for little hybrid paths with a special step
  in their first components is $x (B(x) - S(x)) E(x)$.

  On the other hand, if $P$ doesn't have a special step, then it is a big
  \Sch path and $Q$ can be any little hybrid path.
  The generating function for little hybrid paths with no special step
  in their first components is therefore $x S(x) L(x)$.

  Every nonempty little hybrid path can be uniquely decomposed in this way
  and falls into exactly one of the above categories, so adding in the
  empty path we have
  \begin{displaymath}
    L(x) = 1 + x (B(x) - S(x)) E(x) + x S(x) L(x),
  \end{displaymath}
  or, solving for $L$,
  \begin{equation}
    L(x) = \frac{1 + x ( B(x) - S(x)) E(x)}{1 - x S(x)}.
    \label{eq:L-gf-in-terms-of-B}
  \end{equation}

  The decomposition for big hybrid paths is slightly more involved.
  Given a big hybrid path with an upstep, let $s$ be the first downstep
  to return to the $x$-axis and decompose the path as in
  \autoref{fig:big-hybrid-decomp}.
  In any such big hybrid path, either there is or is not a special
  step before $Q$.

  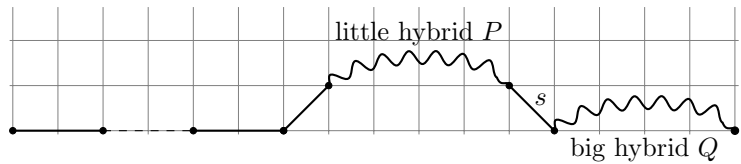
\begin{figure}[h]
      \centering
      \begin{tikzpicture}[scale=3/5]

\begin{scope}[pathfragment]

\startpathat{0}{0} \horiz

\draw[dashed] (2,0) -- (4,0) ;

\startpathat{4}{0} \horiz \up

\draw[pathstep, decorate, decoration={snake}]
    (7,1) .. controls +(45:1) and +(135:1) .. +(4,0) node[midway,above,
    inner sep = 2mm]
    {little hybrid $P$};

\startpathat{11}{1} \down

\node at (11.7,0.7) {$s$};

\node[vertex] at (16,0) {};

\end{scope}


\draw[pathstep, decorate, decoration={snake}]
    (12,0) .. controls +(45:1) and +(135:1) .. +(4,0) node[midway,below,
    inner sep = 4mm]
    {big hybrid $Q$};

\end{tikzpicture}

      \caption{A decomposition of a big hybrid path with an upstep. The
        path begins with a possibly empty sequence of horizontal steps.
        The step $s$ is the first downstep to return to the $x$-axis;
        it may or may not be a special step.}
      \label{fig:big-hybrid-decomp}
  \end{figure}

  Assume that there is a special step before $Q$, and that $s$ is special.
  In that case, $P$ can be any little hybrid path, because if $s$ is
  special, we know that no horizontal steps can appear in $P$ at height
  $1$---such a step has $s$ as its matching downstep and can only
  be created with a slide operation that would make $s$ an ordinary
  step.
  Since there are no horizontal steps at height $1$, $P$ can be any
  little hybrid path, and since there are special steps preceding $Q$, it
  cannot have any horizontal steps, and hence is an odd-special Dyck path.
  The generating function for the possibly empty sequence of horizontal
  steps at the beginning is $1/(1-x)$, so the generating function for
  paths of this type is $x L(x) O(x) / (1 - x)$.

  If a special step appears before $Q$ and $s$ is ordinary, then $P$
  must have a special step and cannot have a horizontal step at height
  ~$1$.
  Assume there is such a horizontal step ~$h$.
  The step $h$ is at odd height and its matching step is $s$, so the
  only way $h$ could be created is by a slide operation that converts $s$
  from a special to an ordinary step, but since we process special steps
  left to right, $s$ would be converted from special to ordinary only if
  there were no special steps in $P$, which is a contradiction.
  This means that $P$ can be any little hybrid path with a special step;
  the generating function for such paths is $L(x) - s(x)$.
  The subpath $Q$ can therefore be any odd-special Dyck path, so the
  generating function for all such big hybrid paths is $x (L(x) - s(x))
  O(x) / (1-x)$.

  Finally, if there is no special step preceding $Q$, then $P$ can be any
  big \Sch path, and $Q$ can be any big hybrid path.
  The generating function for such big hybrid paths is $x S(x) B(x) / (1
  -x)$.

  Every big hybrid path with an upstep falls into exactly one of the
  categories above, so, including paths that consist only of a sequence
  of horizontal steps on the axis, we have
  \begin{displaymath}
    B(x) = \frac{1}{1-x} + \frac{x L(x) O(x)}{1-x} + \frac{x (L(x) - s(x))
      O(x)}{1-x} + \frac{x S(x) B(x)}{1 -x},
  \end{displaymath}
  or, solving for $B$,
  \begin{equation}
   B(x) = \frac{1 + x L(x) O(x) + x (L(x) - s(x))}{1 - x - x S(x)}.
  \label{eq:B-gf-in-terms-of-L}
  \end{equation}

  Solving the system of equations \eqref{eq:L-gf-in-terms-of-B} and
  \eqref{eq:B-gf-in-terms-of-L} and using the fact that
  \begin{gather*}
    s(x) = E(x) = \frac{2}{1+x+R} = \frac{1+x-R}{4x} \quad \text{and}\\
    S(x) = O(x) = \frac{4}{1+x+R} - 1 = \frac{1+x-R}{2x} - 1,
  \end{gather*}
  we obtain the desired expressions for $L(x)$ and $B(x)$.
\end{proof}

It may seem that the generating functions $L$ and $B$ were described in
equations \eqref{eq:little-hybrid-gf-R} and \eqref{eq:big-hybrid-gf-R}
in an unusual way, but the expressions show that $L$ and $B$ are in some
sense built out of familiar generating functions for paths:
\begin{displaymath}
  \frac{2(R+x)}{R (R+x+1)} = 1 + 2x + 7x^{2} + 30x^{3} + 141x^{4}
  + \cdots
\end{displaymath}
is the generating function for sequence A116363, which counts dot
products of rows of Pascal's and Catalan's triangle, and of course
$2/(R+x+1)$ is the generating function for the little \Sch numbers.
Also appearing in both $L$ and $B$ is
\begin{displaymath}
  \frac{R+1-x}{2} = 1 - x - x S(x),
\end{displaymath}
a minor modification of the generating function for the big \Sch
numbers.
In $B$, we see that we have exactly $-x S(x)$; the remaining terms in
$B$ are
\begin{displaymath}
  -\frac{R(R+x)}{2} + \frac{3}{2} - \frac{3x}{2}
\end{displaymath}
which is $1 + x + x^{2} S(x)$.

\section{The bijections $E$ and $e$ and 231-avoiding permutations}
\label{sec:hybrid-paths-231}

While using the bijections $E$ or $e$, one can keep track of the order
in which horizontal steps are added and thereby associate a permutation
to an even- or odd-special Dyck path.
For example, the paths in \autoref{fig:E-forward-backward} correspond to
the permutations $321$, $123$, and $132$; the path in
\autoref{fig:e-bijection-ex} corresponds to $12387465$.
In this section we will see that every permutation so obtained must avoid
the pattern $231$. This is very interesting, since $231$-avoiding
permutations are counted by the Catalan numbers and hence are in
bijection with Dyck paths; see Mansour et al. \cite[\S
3.1]{mansour.deng.ea:dyck}.

A permutation $\pi = \pi_{1}\pi_{2}\cdots \pi_{n}$ written in one-line
notation \textsl{contains} a pattern $\sigma$ (another permutation) if
there is some subset of the $\pi_{i}$'s that are order-isomorphic to
$\sigma$.
A permutation \textsl{avoids} a pattern if it does not contain it.
The permutation $12584367$ contains the pattern $231$ because the subset
$583$ is order-isomorphic to $231$, and avoids the pattern $3124$.
The notation $S_{n}(231)$ refers to the set of $231$-avoiding
permutations of $[n]$.

We start with a lemma that tells us exactly when the horizontal steps
created by two special steps are added out of order---that is, when
two special steps create the pattern $21$.
Given a special step $s$, let $h(s)$ refer to the horizontal step
created when $s$ is turned into an ordinary step.

\begin{lemma}
  \label{thm:21-order-adding-sp-steps}
  Given two special steps $a$ and $b$ in a hybrid path with $a$ to the
  left of $b$, $h(b)$ is created to the left of $h(a)$ if and only if
  $b$ is preceded by a downstep $d$ and the upstep matching $d$ is to
  the left of $a$.
\end{lemma}

\begin{proof}
  We only need to examine three possibilities: $b$ is preceded by an
  upstep, $b$ is preceded by a downstep whose matching upstep is to
  the right of $a$, and $b$ is preceded by a downstep whose matching
  upstep is to the left of $a$.
  To work through those three cases, we need to use the fact that for
  any special step $s$, $h(s)$ is created to the left of $s$ and to
  the right of the downstep matching $s$.
  (When doing a flatten, ``to the left'' and ``to the right'' are weak
  inequalities, since the horizontal step will be created in the same
  position as those steps.)

  In the first case, if $b$ is preceded by an upstep, then $h(b)$ will
  clearly be to the left of $h(a)$, since $h(b)$ will be created at
  the position of $b$, which is to the right of $a$.

  If $b$ is preceded by a downstep whose matching upstep is to the
  right of $a$, let $u$ be that matching upstep.
  Since $h(b)$ will be created at $u$ and the following step, $h(b)$
  is to the right of $a$ and hence to the right of $h(a)$.

  Finally, if $b$ is preceded by a downstep whose matching upstep is
  to the \emph{left} of $a$, let $u$ be that matching upstep.
  See \autoref{fig:e-bijection-slide}; $b$ would be the special step
  pictured in that figure, $u$ would be the upstep, and $a$ would be
  somewhere in the subpath $P$, and since the upstep matching $a$ is
  also in that subpath, $h(b)$ will be created to the left of $h(a)$.
\end{proof}

With that result, we can easily prove the following theorem.

\begin{theorem}
  \label{thm:231-avoiding}
  The permutation corresponding to the order in which horizontal steps
  are added while transforming an even- or odd-special Dyck path into a
  small or large \Sch path avoids the pattern~$231$.
\end{theorem}

\begin{proof}
  Consider any three special steps $a$, $b$, and $c$ in a hybrid path,
  appearing in that order left to right.
  If these three steps cause the corresponding permutation to contain
  $231$, then we must have $h(b)$, $h(c)$, and $h(a)$ in that order.
  The bijections $E$ and $e$ process special steps left to right, so we
  first create $h(a)$, and then create $h(b)$ to the left of that.
  Now we must have $h(c)$ created to the right of $h(b)$, which means
  by \nonthmref{thm:21-order-adding-sp-steps} either $c$ is preceded by an
  upstep, or is preceded by a downstep whose matching upstep is to the
  right of $b$, but both of those possibilities cause $h(c)$ to be to
  the right of $h(a)$, which is a contradiction.
\end{proof}

The permutations produced are therefore a subset of $231$-avoiding
permutations; next we will see that every such permutation can be
obtained from some even- or odd-special Dyck path.

\begin{theorem}
  \label{thm:hybrid-paths-all-231-avoiding-perms-possible}
  Every $231$-avoiding permutation can be obtained from some
  odd-special Dyck path.
\end{theorem}

\begin{proof}
  Given a $231$-avoiding permutation $\pi = \pi_{1}\cdots \pi_{n}$, we
  use the following recursive procedure to construct an odd-special Dyck
  path that, when using $E_{\infty}$, will create horizontal steps in the
  order specified by $\pi$.

  To begin with, the empty permutation corresponds to the empty path.
  Given a nonempty permutation $\pi$ of length $n$, if $\pi_{n} = n$,
  find the path corresponding to $\pi_{1}\cdots\pi_{n-1}$ and
  append the steps Ud to that path.
  (Recall that U refers to an upstep, and D and d refer to ordinary
  and special downsteps, respectively.)
  We will call this an \textsl{append} operation; it corresponds to the
  flatten operation.
  If a path $P$ corresponds to the permutation $\pi_{1}\pi_{2}\dots
  \pi_{n-1}$, then the path obtained by appending Ud to $P$ will
  correspond to $\pi_{1}\pi_{2}\dots \pi_{n-1} n$.

  If the permutation does not end with $n$, we will need the
  \textsl{lift} operation, which is defined as follows: given an
  odd-special Dyck path of the kind produced by this procedure, find
  an upstep leaving from the $x$-axis and let $Q$ be the subpath
  consisting of that upstep and everything following it.
  The lift operation, illustrated in \autoref{fig:lift}, replaces $Q$
  with a path consisting of two upsteps, then $Q$, then an ordinary
  downstep and a special downstep.
  The lift operation corresponds to the slide operation of $E$ and
  $e$.

  \begin{figure}[h]
      \centering
      \begin{tikzpicture}[scale=.66]

\begin{scope}[pathfragment]
  \draw[pathstep, decorate, decoration={snake}]
    (0,0) .. controls +(45:1) and +(135:1) .. +(3,0) node[midway,above,
    inner sep = 2mm] {$P$};

  \filldraw (3,0) circle (2pt);

  \draw[pathstep, decorate, decoration={snake}]
    (3,0) .. controls +(45:1) and +(135:1) .. +(3,0) node[midway,above,
    inner sep = 2mm] {$Q$};

\end{scope}

\draw[->] (6.5,2/3) -- (7.5, 2/3);

\begin{scope}[pathfragment, xshift=8cm]

  \draw[pathstep, decorate, decoration={snake}]
    (0,0) .. controls +(45:1) and +(135:1) .. +(3,0) node[midway,above,
    inner sep = 2mm] {$P$};

  \startpathat{3}{0} \up \up

  \draw[pathstep, decorate, decoration={snake}]
    (5,2) .. controls +(45:1) and +(135:1) .. +(3,0) node[midway,above,
    inner sep = 2mm] {$Q$};

  \startpathat{8}{2} \down \spdown

\end{scope}

\end{tikzpicture}

      \caption{The lift operation. If the resulting path has $n$ special
        steps, in the corresponding permutation all numbers from $P$
        will precede $n$, and all numbers from $Q$ will follow $n$.}
      \label{fig:lift}
  \end{figure}
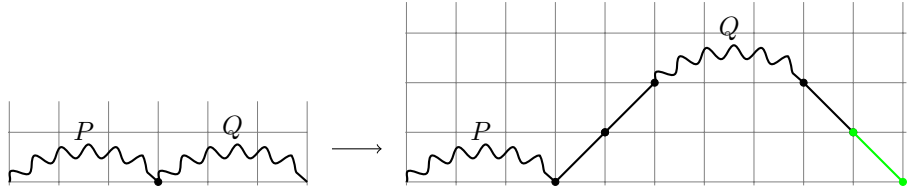

  \newcommand{\corrto}{\leftrightarrow}

  For ease of description, define \textsl{good insertion} to be the
  operation of inserting $n$ into a $231$-avoiding permutation of
  $[n-1]$ anywhere except at the end so that the resulting
  permutation is also $231$-avoiding.
  We write $P \corrto \pi$ if a path $P$, constructed using the append
  and lift operations, corresponds to the permutation $\pi$.

  Assume $P \corrto \pi$, where $\pi \in S_{n}(231)$.
  Say that we obtain $\pi'$ by good insertion of $n+1$ after the $k$th
  entry of $\pi$.
  We need to show first that we can do a lift operation following the
  $k$th special step of $P$, and second that the corresponding path $P'$
  corresponds to $\pi'$.

  The path Ud corresponds to the permutation $1$, and both claims are
  true for that path-permutation pair.
  To prove the two claims in general, we need two propositions:

  \begin{proposition}
    If $\pi \in S_{n}(231)$, good insertion can be done after the $k$th
    entry of $\pi$ if and only if the first $k$ entries of $\pi$ form a
    permutation of $[k]$.
  \end{proposition}

  \begin{proposition}
    Assume that $P \corrto \pi$.
    The $k$th special step of $P$ ends on the $x$-axis if and only if
    the first $k$ entries of $\pi$ form a permutation of $[k]$.
  \end{proposition}

  The proof of the first is elementary and left to the reader.
  As for the second, let $s$ be the $k$th special step of $P$.
  Assume $s$ ends on the $x$-axis.
  Then by \nonthmref{thm:21-order-adding-sp-steps}, the horizontal step
  for every special step to the right of $s$ will be created to the right
  of $s$, which means the first $k$ entries of $\pi$ form a permutation of
  $[k]$.
  On the other hand, if $s$ does not end on the $x$-axis, then because
  of the definition of the lift operation, there must be a special step to
  the right of $s$ that is immediately preceded by a downstep $d$, with
  the upstep matching $d$ to the left of $s$.
  Therefore, by \nonthmref{thm:21-order-adding-sp-steps}, there will be a
  number bigger than $k$ among the first $k$ entries of $\pi$, so
  $\pi_{1}\cdots \pi_{k}$ will not form a permutation of $[k]$.

  The first claim above is now clear.
  The second claim is also easy to see: say $\pi = \pi_{1}\cdots
  \pi_{n}$ and $\pi'$ is obtained by good insertion after $\pi_{k}$.
  If one lifts $P$ after the $k$th special step, the resulting path $P'$
  will correspond to the permutation $\pi'$ because, following the lift
  operation, the special steps corresponding to $\pi_{1}\cdots\pi_{k}$
  and to $\pi_{k+1}\cdots \pi_{n}$ will be turned into horizontal steps
  in the same order (they are either unchanged or simply raised by $2$
  units), and the final step of $P'$ is a special downstep that will be
  turned into a horizontal step that follows every horizontal step
  created by the first $k$ special steps of $P'$ and precedes every
  horizontal step created by the $(k+1)$st to $n$th special steps.
\end{proof}

\begin{figure}[h]
  \centering
  \begin{tikzpicture}[scale=.29]

\begin{scope}[pathfragment]
  \startpathat{0}{0} \up \spdown \up \up \up \up \up \spdown \down
  \spdown \up \spdown \down \spdown
\end{scope}

\node at (7, -1) {$15324$};
\draw[<-] (14.5,1.5) -- (15.5, 1.5);

\begin{scope}[pathfragment, xshift=16cm]
  \startpathat{0}{0} \up \spdown \up \up \up \spdown \down \spdown \up
  \spdown
\end{scope}

\node at (21, -1) {$1324$};
\draw[<-] (26.5,1.5) -- (27.5, 1.5);

\begin{scope}[pathfragment, xshift=28cm]
  \startpathat{0}{0} \up \spdown \up \up \up \spdown \down \spdown
\end{scope}

\node at (32, -1) {$132$};
\draw[<-] (36.5, 0.5) -- (37.5, 0.5);

\begin{scope}[pathfragment, xshift=38cm]
  \startpathat{0}{0} \up \spdown \up \spdown
\end{scope}

\node at (40, -1) {$12$};
\draw[<-] (42.5, 0.5) -- (43.5, 0.5);

\begin{scope}[pathfragment, xshift=44cm]
  \startpathat{0}{0} \up \spdown
\end{scope}

\node at (45, -1) {$1$};
\draw[<-] (46.5, 0.5) -- (47.5, 0.5);

\begin{scope}[pathfragment, xshift=49cm]
  \fill (0,0) circle (3pt);
  \draw (-1/2, 1) circle (0pt);
  \draw (1/2, 1) circle (0pt);
\end{scope}

\node at (49, -1) {$\varnothing$};

\end{tikzpicture}

  \caption{An example of the recursive procedure to build an OSDP
    corresponding to a $231$-avoiding permutation. Here we see how the
    path for $15324$ is built up from the empty path. Below each path
    is the corresponding permutation.}
  \label{fig:permutation-recursion}
\end{figure}
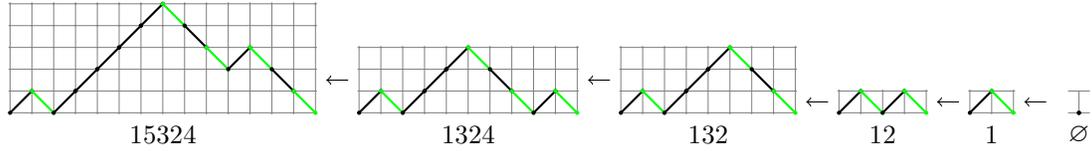

\autoref{fig:permutation-recursion} shows an example of this procedure.
One can also obtain every $231$-avoiding permutation with even-special
Dyck paths simply by taking a path produced by this procedure and
sandwiching it between an upstep and an ordinary downstep.
That produces an even-special Dyck path that clearly corresponds to
the same permutation.

We close with an interesting conjecture.
The paths produced for $231$-avoiding permutations of $[n]$ are not
all the same length; the lengths range from $2n$ for the path
corresponding to $123\cdots n$ to $4n-2$ for the path corresponding to
$n (n-1)\cdots 2 1$.
An obvious question to ask is: how are the lengths distributed?
In other words, find the coefficients of the polynomial
\begin{displaymath}
  \sum_{\pi \in S_{n}(231)} q^{\pathlen(\pi)},
\end{displaymath}
where $\pathlen(\pi)$ is the length of the path corresponding to $\pi$
using the construction above.
We can also sum that expression over all $n$, since for a given length
there can be only finitely many permutations that correspond to a path
of that length, and ask what generating function we get.

The lengths appear to have the Narayana distribution; see sequence
A1263 and Sulanke~\cite{sulanke:narayana}. The table below shows the
polynomials for some small values of $n$.

\begin{table}[h]
  \centering
  \begin{tabular}{ll|ll}
    $n$ & distribution of lengths & $n$ & distribution of lengths\\
    \hline
    $1$ & $q^{2}$  & 4 & $q^{14} + 6q^{12} + 6q^{10} + q^{8}$ \\
    $2$ & $q^{6} + q^{4}$ & 5 & $q^{18} + 10q^{16} + 20q^{14} + 10q^{12} + q^{10}$ \\
    $3$ & $q^{10} + 3 q^{8} + q^{6}$ & 6 &  $q^{22} + 15q^{20} + 50q^{18} + 50q^{16} + 15q^{14} + q^{12}$ \\
  \end{tabular}
    \caption{The lengths of the OSDPs corresponding to $231$-avoiding
      permutations of $[n]$ appear to be Narayana-distributed.}
  \label{tab:lengths-narayana-dist}
\end{table}

\begin{conjecture}
  The lengths of the paths corresponding to $231$-avoiding permutations
  of $[n]$ using the above construction have the Narayana distribution;
  that is,
  \begin{equation}
      \sum_{\pi \in S_{n}(231)} q^{\pathlen(\pi)} =
      q^{2n} \sum_{k \ge 0} N(n, k) q^{2k},
  \label{eq:conj-dist}
  \end{equation}
  where $N(n,k)$ is a Narayana number. We also have
  \begin{equation}
    \sum_{n \ge 0} \  \sum_{\pi \in S_{n}(231)} q^{\pathlen(\pi)} =
    \frac{1-q^{2}+q^{4}-\sqrt{1-2q^{2}-q^{4}-2q^{6}+q^{8}}}{2q^{4}}.
  \label{eq:conj-gf}
  \end{equation}
\end{conjecture}

The right-hand side of equation~\eqref{eq:conj-gf} is the generating
function for generalized Catalan numbers described by Stein and
Waterman~\cite{stein.waterman:new} (see the $m=1$ column of Table 1),
and by Vauchassade de Chaumont and Viennot
\cite{chaumont.viennot:polynomes}. Those numbers are sequence A4148,
and count secondary structures of RNA molecules according to the number
of bases.

\section{Included Sage code}
\label{sec:included-sage-code}

\newcommand{\myurl}[1]{\href{http://#1}{\texttt{#1}}}

This preprint includes Sage code
(see \myurl{sagemath.org}) for working with the various
paths, bijections, matchings, and permutations described here.
The code is included with the preprint source; visit
\myurl{arxiv.org/abs/1006.1959}, click ``Other formats'', then
``Download source'', and look for the file
\texttt{code\_for\_dyck\_schroeder.sage}.

\bibliographystyle{amsplainurl} \bibliography{dyck-schroeder}

\end{document}